\documentclass{article}
\usepackage{amsmath,amsthm}
\usepackage{amsfonts}
\newtheorem{theorem}{Theorem}

\newtheorem{lemma}{Lemma}
\newtheorem{proposition}[theorem]{Proposition}

\begin{document}
\bibliographystyle{plain}
\begin{center}
{\LARGE \bf A stochastic flow  arising in the study of local times}\\
\vspace{.1in}
{\large JON WARREN}\\
\vspace{.1in} Department of Statistics, University of Warwick, Coventry CV4 7AL, UK.
\end{center}
\begin{abstract}
A stochastic flow of homeomorphisms of ${\mathbb R}$ previously
studied by Bass and  Burdzy  \cite{bb} and Hu and Warren \cite{hw} is
shown to arise in the study of the local times of Brownian
motion. This leads to  a new proof of the Ray-Knight theorems for the flow via the classical Ray-Knight theorems for Brownian motion.
\end{abstract}
\section{Introduction}

Let $\beta_1$ and $\beta_2$ be fixed real constants. Suppose $B$ is a
Brownian motion on  the real line issuing from $0$. Associated with
the equation 
\begin{equation}
\label{flow}
X_t(x)= x+ B_t + \beta_1 \int_0^t ds {\mathbf 1}{(X_s(x) \leq 0)} + \beta_2
\int_0^t ds {\mathbf 1}{(X_s (x)>0)},
\end{equation}
is a stochastic flow. There exists a random flow of homeomorphisms of
the real line $(X_t;t \geq 0)$ so that  $t \mapsto X_t(x)$ satisfies
the above equation for all $x \in {\mathbb R}$.
For each $t$ the map $x \mapsto X_t(x)$ is increasing and
differentiable. We denote the derivative by $DX_t(x)$.
It has been  shown by Bass and Burdzy, \cite{bb}, and Hu and Warren, \cite{hw},  that for certain stopping times $T$ the
process $x \mapsto DX_T(x)$ is a diffusion process. Such a  result  was
called a Ray-Knight theorem for the flow.

Let $(W_t; t \geq 0)$ be a real-valued Brownian motion starting from
$W_0=\xi>0$. Let $T_0$ be the stopping time
\begin{equation}
T_0= \inf\{t \geq 0: W_t=0\}.
\end{equation}
Let $(l(t,x); t\geq 0, x \in {\mathbb R})$ be the family of
bi-continuous local times of $W$.  The celebrated Ray-Knight theorem
states that the process $ x \mapsto l(T_0,x)$ is a diffusion.
 The problem studied by Warren and Yor \cite{wy} (actually a slight
 variant of it)  was that of obtaining a
 description of $(W_t; 0 \leq t \leq T_0)$ after conditioning on the
 values of the entire family of random variables $(l(T_0,x); x \in
 {\mathbb R})$. A related problem was treated by Aldous, \cite{aldous}.
The method of \cite{wy} was to write $W$ as a random 
transformation of a  process $\hat{W}$, christened the Burglar process, which is 
independent of $(l(T_0,x); x \in
 {\mathbb R})$. It was   shown, using a skew product linking squared Bessel and Jacobi diffusions,   that Ray-Knight type theorems hold
 for  $\hat{W}$ itself.
The principal result of this paper is the following theorem which
establishes a close connection between the this problem and a  stochastic
flow of the Bass-Burdzy type. 
Let $l(t,x)$ be the local times of a real-valued Brownian motion $W$
as above. Define 
\begin{equation}
\lambda(t,x)= l(T_0,x)- l(t,x), \qquad \text{ for } x \in{\mathbb R}, \qquad  t
\in [0, T_0].
\end{equation}
Let $M= \sup\{ W_t: t \in [0,T_0]\}$ and $T_M= \inf\{ t: W_t=M\}$. Now
define a family of random changes of scale via:
\begin{equation}
\phi_t(x)= \int_x^{W_t} \frac{dz}{\lambda(t,z)} \qquad \text{ for } x \in (0,M), \qquad t
\in [0, T_M).
\end{equation}
Each map $\phi_t$ is a homeomorphism between $(0,M)$ and ${\mathbb R}$ and
satisfies $\phi_t(W_t)=0$.
Next define a random clock via:
\begin{equation}
A_t= \int_0^t \frac{ds}{\lambda(s,W_s)^2} \qquad \text{ for } 0 \leq t<T_M.
\end{equation}
We will see below that $A_{{T_M}-}= \infty$, almost surely.
\begin{theorem}
\label{main}
The family of maps $(X_u; 0 \leq u <\infty)$ defined via
\[
X_{A_t}(x)= \phi_t \circ \phi_0^{-1}(x) \qquad \text{ for } x \in
{\mathbb R},\qquad  0 \leq t<{T_M},
\]
 is distributed as the
Bass-Burdzy flow associated with \eqref{flow} with $\beta_1=0$ and
$\beta_2=1$. Moreover it is independent of  $(l(T_0,x); x \in
 {\mathbb R})$.
\end{theorem}

\section{From the Burglar to a flow}

In this section  we prove Theorem \ref{main}. According to this result, the flow $(X_u; 0 \leq u <\infty)$ constructed from the Brownian motion $W$ 
 is independent of the local times 
 $(l(T_0,x); x \in
 {\mathbb R})$. Because of this, Theorem \ref{main} is intimately related to the problem considered in \cite{wy} of describing the process $W$ conditional 
 on knowing these local times.  One standard approach for studying  this type of problem is  the method of enlargement of filtrations. 
 We assume that the Brownian motion $\bigl(W_t; t\geq 0\bigr)$ is carried by a filtered probability 
 space $\bigl(\Omega,({\cal F}_t)_{t\geq 0},{\mathbb P} \bigr)$.  Let ${\cal L}= \sigma\bigl( l(T_0, x); x \in{\mathbb R}\bigr)$ and introduce 
 the enlarged filtration $\bigl( {\cal F}^\ast_t; t\geq 0\bigr)$ where ${\cal F}^\ast_t= {\cal F}_t \vee {\cal L}$.   Typically we would seek  to show that a martingale in the original 
 filtration is a semimartingale with respect to the enlarged filtration, and to obtain some expression for its canonical decomposition as such. However in this case there
 is no reason to believe that $W$  is a semimartingale with respect to ${\cal F}^\ast_t$;  Jacod's   absolute continuity criteria, \cite{j}, does not hold.  
 Thus another approach is required to generate examples of martingales and semi-martingales relative to ${\cal F}^\ast_t$.

  For $x\in {\mathbb R}$
  we denote by $T_x$ the first time $W$ reaches  the level $x$. We adopt the usual convention that $\exp(-\infty)=0$.
\begin{proposition}
\label{mart}
 Let $\kappa> 2$ be a  constant, and let $a$,$\xi$ and $c$  be constants satisfying $W_0=\xi$ and  $0<a<\xi<c$.
The  process 
\[
M^\ast_t=   \exp
\left\{\frac{-\kappa^2}{2}\int_a^{W_t} \frac{dy}{\lambda(t,y)} - \frac{\kappa(\kappa-2)}{2}\int_{W_t}^c \frac{dy}{\lambda(t,y)} \right\} \qquad t\leq T_a \wedge T_c, 
\]
$M^\ast_t= M^\ast_{T_a\wedge T_c}$ for $t>T_a\wedge T_c$,  is an ${\cal F}^\ast_t$-martingale.
\end{proposition}

Admitting this proposition for the moment, let us see how it leads to a proof the main result. First note:

\begin{lemma}
\label{occ}
\[
\int_a^c \frac{dz}{\lambda(t,z)} -\int_a^c \frac{dz}{\lambda(0,z)}= \int_0^t {\mathbf 1}{\bigl( a\leq W_s \leq c\bigr)}dA_s 
\]
almost surely on the event $\{ c<M, t< T_M\}$.
\end{lemma}
\begin{proof} Notice  that on the event $\{c<M, t<T_M\}$ we have, almost surely, $\lambda(s,z)$ being bounded away 
from $0$ for all $z \in [a,c]$ and $s\in [0,t]$.
Then the  result is a consequence of applying the  occupation time formula, see Exercise (1.15), Chapter VI of Revuz and Yor, \cite{ry}, to obtain:
\[
\int_a^c dz\int_0^t \frac{d_sl(s,z)}{\lambda(s,z)^2}=\int_0^t {\mathbf 1}{\bigl( a\leq W_s \leq c\bigr)}\frac{ds} {\lambda(s,W_s)^2}.
\]
\end{proof}

We can use this lemma to verify that $A_{T_M-}=\infty$ almost surely. For we see that on $\{c<M\}$,
\begin{equation}
\label{infinity}
A_{T_M-}\geq 
\int_a^c \frac{dz}{\lambda(T_M,z)} -\int_a^c \frac{dz}{\lambda(0,z)}.
\end{equation}
The distributions of $\lambda(T_M, \cdot)$ and $\lambda(0,\cdot)$ can be described  using the classical Ray-Knight theorems for Brownian motion, 
as in Lemma \ref{classic}, and then it follows from asymptotics of squared Bessel processes, to be found in for instance in  \cite{ry}, that almost surely,
\begin{equation}
\lim_{\epsilon \downarrow 0} \frac{1}{\log(1/\epsilon)} \int_a^{M-\epsilon} \frac{dz}{\lambda(T_M,z)}= \infty \qquad \text{ and } 
\qquad \lim_{\epsilon \downarrow 0} \frac{1}{\log(1/\epsilon)} \int_a^{M-\epsilon} \frac{dz}{\lambda(0,z)}= 1/2.
\end{equation}

We define a process $\bigl(B_t; t\geq 0\bigr)$ on the sample space $\Omega$ which carries $W$ via
\begin{equation}
\label{bdef}
B_{A_t}= \int_{W_0}^c \frac{dz}{\lambda(0,z)}-\int_{W_t}^c \frac{dz}{\lambda(t,z)} \qquad \text{ on }\{ t< T_c\} \cap \{c<M\}.
\end{equation}
Note that this equation defines $B$ consistently as $c$ varies, and 
that $B$ is defined for arbitrarily large times by letting $c \uparrow M$ so that $T_c \uparrow T_M$ and $A_{T_c} \uparrow \infty$.
 
Let $X$ be defined as in Theorem \ref{main}. Suppose that $0<x<M$ and $t<T_M$.  Then  for $c\in (x,M)$ satisfying $t<T_c$,
\begin{equation}
\label{key}
\begin{split}
X_{A_t}\bigl(\phi_0(x)\bigr)-\phi_0(x) &= \int_x^{W_t} \frac{dz}{\lambda(t,z)}- \int_x^{W_0}\frac{dz}{\lambda(0,z)} \\
&=  \int_x^{c} \frac{dz}{\lambda(t,z)}- \int_{W_t}^c \frac{dz}{\lambda(t,z)}-
\int_x^{c}\frac{dz}{\lambda(0,z)}+\int_{W_0}^c\frac{dz}{\lambda(0,z)} \\
&=B_{A_t} +\int_x^{c} \frac{dz}{\lambda(t,z)}-\int_x^{c}\frac{dz}{\lambda(0,z)}
\\
&=B_{A_t} +\int_0^t {\mathbf 1}\bigl( X_{A_s}(\phi_0(x)) \geq 0\bigr) dA_s.
 \end{split}
\end{equation}
Here on the last line we have used Lemma \ref{occ} together with the fact that for $s<T_c$,
\begin{equation}
 {\mathbf 1}{\bigl( x\leq W_s \leq c\bigr)}= {\mathbf 1}{\bigl( x\leq W_s \bigr)}=  {\mathbf 1}{\bigl( X_{A_s}(\phi_0(x)) \geq 0  \bigr)}.
\end{equation}
Time-changing equation \eqref{key},  and putting $y=\phi_0(x)$, we see that $X_t(y)$  solves  equation \eqref{flow}  driven by the process $B$.

 We must next show that $B$ is a Brownian motion with respect to a suitable filtration. Using Lemma \ref{occ} we can re-write the martingale $M^\ast_t$ as,
\begin{equation}
\frac{M^\ast_t}{M^\ast_0}= \exp\bigl\{-\kappa B_{A_t}-\tfrac{1}{2}\kappa^2A_t \bigr\} \qquad \text{ on } \{ t<T_a \wedge T_c\} \cap \{c<M\}.
\end{equation}
By letting $a$ decrease to  $0$ and $ c$ increase to $M$ we see that $\exp\bigl\{-\kappa B_{A_t}-\tfrac{1}{2}\kappa^2A_t \bigr\}$ 
is  ${\cal F}^\ast_t$-local martingale on the stochastic interval $[0,T_M)$. Note here that $T_M$ is an ${\cal F}^\ast_t$-stopping time.
Introduce that filtration $\bigl({\cal G}_u; u\geq 0\bigr)$ via
\begin{equation}
{\cal G}_u = {\cal F}^\ast_{\alpha_u}, \qquad u\in [0,\infty),
\end{equation}
where $(\alpha_u;u \geq 0)$ is the inverse of $(A_t;t\geq 0)$.
The local martingale property is preserved by the time-change and we are able to conclude that $B$ is a Brownian motion by virtue of the next lemma. 

\begin{lemma} If $(B_t; t\geq 0)$ is a continuous ${\cal G}_t$-adapted process such that 
\[
\exp\bigl\{-\kappa B_t-\tfrac{1}{2}\kappa^2 t\bigr\} \text{ is a } {\cal G}_t\text{-local martingale}
\]
for at least two distinct and non-zero values of $\kappa$, then $(B_t; t\geq 0)$ is a  ${\cal G}_t$-Brownian motion.
\end{lemma}
\begin{proof}
By considering logarithms it is clear that $B_t$ is a ${\cal G}_t$ semimartingale. Let its canonical decomposition be $B_t=B_0+N_t+A_t$ where $A_t$ is a finite variation process and $N_t$ is a ${\cal G}_t$-local martingale. Applying It\^{o}'s formula to $Z^{(\kappa)}_t=\exp\{-\kappa B_t -\tfrac{1}{2}\kappa^2 t\}$, we obtain 
\[
Z^{(\kappa)}_t= Z^{(\kappa)}_0- \kappa\int_0^t Z^{(\kappa)}_s\; dB_s + \kappa^2/2 \int_0^t Z^{(\kappa)}_s\;(d[N]_s-ds).
\]
Using this we see that the finite variation part of the local martingale $\tfrac{1}{\kappa}\int_0^t dZ^{(\kappa)}_s/Z^{(\kappa)}_s$, which must be identically zero, is given by 
\[
-A_t+\tfrac{\kappa}{2}(d[N]_s-ds).
\]
From this we conclude that $A_t$ is identically zero, that is to say $B_t$ is a local martingale itself, and that $[B]_t=[N]_t=t$. Thus $B$ is a ${\cal G}_t$-Brownian motion by virtue of It\^{o}'s characterization.
\end{proof}

Finally to complete the proof of Theorem \ref{main} we must note that the Brownian motion $B$ is independent of the $\sigma$-algebra ${\cal G}_0$ which contains ${\cal L}$, and that since  the stochastic differential equation  \eqref{main} is exact by Zvonkin's observation,  the flow $X$ is measurable with respect to $B$, and hence also independent of ${\cal L}$.

We turn now to proving Proposition \ref{mart}.  The principle tool used in constructing ${\cal F}^\ast_t$-martingales is a certain  Markov property.
Assume that   the sample space $\Omega$ is the canonical space ${\cal C}\bigl([0,\infty), {\mathbb R}\bigr)$ and $\bigl( W_t; t\geq 0\bigr)$ realized as the co-ordinate process on $\Omega$. Then we can introduce the family $\bigl( \theta_t ; t \geq 0 \bigr)$ of shift operators:
\[
W_t \circ \theta_s =W_{t+s}, \qquad \qquad s,t \geq 0.
\]
For $x \in {\mathbb R}$ let  $ {\mathbb P}^x$ be the  probability measures on $\Omega$ under which  the co-ordinate process $W$ is a Brownian motion starting from $x$.   Then
let ${\mathbb P}^{x,\ell}$ for $\ell \in L$ denote a  version of the  regular conditional probability for  ${\mathbb P}^x$  given $l(T_0, \cdot)= \ell (\cdot)$, where $L$ is some suitable space of local time profiles.   
\begin{proposition}
The process $\bigl( W_t, \lambda(t, \cdot); 0\leq t <T_0\bigr)$ is Markovian relative to the filtration $\bigl( {\cal F}^\ast_t; t\geq 0 \bigr)$:
\[
{\mathbb E}^x\bigl[ F\circ \theta_t \vert {\cal F}_t^\ast \bigr] = {\mathbb P}^{W_t,\lambda(t,\cdot)}\bigl(F\bigr), \qquad a.s. \text{ on } \{ t<T_0\},
\]
for all $x>0$, and  non-negative, measurable $F$ on $\Omega$.
\end{proposition}
\begin{proof}
Let $\nu$ be a finite measure supported in $(0,\infty)$. 
Observe that, on $\{ t<T_0\}$,
\[
{\cal E}_\nu\bigl(( l(T_0,\cdot)\bigr)= {\cal E}_\nu\bigl((\lambda(t,\cdot)\bigr){\cal E}_\nu\bigl(l(t, \cdot)\bigr),
\]
where  ${\cal E}_\nu(l)$ denotes
\[
 \exp \left\{ -\tfrac{1}{2} \int_0^\infty \nu(dy)l(y) \right\}.
\]
If we define $F^\nu$ to be $F {\cal E}_\nu\bigl(l(T_0, \cdot)\bigr)$ then
\[
F^\nu\circ \theta_t \;{\cal E}_\nu\bigl(l(t, \cdot)\bigr) = F\circ \theta_t\; {\cal E}_\nu\bigl( l(T_0, \cdot)\bigr)  \qquad \text { on } \{t <T_0\}.
\]
Also put $G^\nu= {\mathbb P}^{W_0,l(T_0,\cdot)}(F) {\cal E}_\nu\bigl(l(T_0,\cdot)\bigr)$. Then  ${\mathbb E}^y[G^\nu]= {\mathbb E}^y[F^\nu]$, for all $y>0$, and 
\[
G^\nu \circ \theta_t  = {\mathbb P}^{W_t,\lambda(t,\cdot)}(F) {\cal E}_\nu\bigl( \lambda(t,\cdot)\bigr)  \qquad \text { on } \{t <T_0\}.
\]
Using these various observations and two applications of  the Markov property of Brownian motion we have, for each $A \in {\cal F}_t$
\begin{equation*}
\begin{split}
{\mathbb E}^x\bigl[ F\circ \theta_t {\mathbf 1}_A {\mathbf 1}_{(t<T_0)} {\cal E}_\nu\bigl(l(T_0,\cdot)\bigr)\bigr] &= 
{\mathbb E}^x\bigl[ F^\nu\circ \theta_t {\mathbf 1}_A {\mathbf 1}_{(t<T_0)} {\cal E}_\nu\bigl(l(t,\cdot)\bigr)\bigr]\\
&={\mathbb E}^x\bigl [ {\mathbb E}^{W_t}\bigl[ F^\nu \bigr] {\mathbf 1}_A{\mathbf 1}_{(t<T_0)}  {\cal E}_\nu\bigl(l(t,\cdot)\bigr)\bigr]
\\
&={\mathbb E}^x\bigl [ {\mathbb E}^{W_t}\bigl[ G^\nu \bigr] {\mathbf 1}_A{\mathbf 1}_{(t<T_0)}  {\cal E}_\nu\bigl(l(t,\cdot)\bigr)\bigr]
\\
&= {\mathbb E}^x\bigl [  G^\nu\circ \theta_t {\mathbf 1}_A{\mathbf 1}_{(t<T_0)}  {\cal E}_\nu\bigl(l(t,\cdot)\bigr)\bigr] \\
&={\mathbb E}^x\bigl [ {\mathbb P}^{W_t, \lambda(t,\cdot)}\bigl( F \bigr) {\mathbf 1}_A{\mathbf 1}_{(t<T_0)}  {\cal E}_\nu\bigl(l(T_0,\cdot)\bigr)\bigr].
\end{split}
\end{equation*}
Since, the random variables ${\mathbf 1}_A$ and ${\cal E}_\nu\bigl(l(T_0,\cdot)\bigr)$ generate ${\cal F}^\ast_t$ as $A$ and $\nu$ vary we are done. 
\end{proof}

 Consider $F$ given by
\begin{equation}
F= {\mathbf 1}_{( T_c<T_a)} \exp
\left\{ \frac{-\kappa^2}{2}\int_a^{c} \frac{dy}{\lambda(T_c,y)} \right\}  
+  {\mathbf 1}_{( T_a<T_c)}  \exp\left\{-\frac{\kappa(\kappa-2)}{2}\int_{a}^c \frac{dy}{\lambda(T_a,y)} \right\},  
\end{equation}
where $0<a<\xi<c$, recalling $W_0=\xi$.
Using the Markov property just proved, the martingale $M^\ast_t$ is constructed as
\begin{equation}
\begin{split}
M^\ast_t= {\mathbb E}^\xi \bigl[ F\vert {\cal F}^\ast_t \bigr]
= {\mathbb E}^\xi \bigl[ F \circ \theta_t {\mathbf 1}_{(t <T_a \wedge T_c)}+ F{\mathbf 1}_{(t >T_a \wedge T_c)} \vert {\cal F}^\ast_t \bigr] \\
= {\mathbb P }^{W_t, \lambda(t, \cdot)}\bigl (F\bigr)  {\mathbf 1}_{(t <T_a \wedge T_c)}+F{\mathbf 1}_{(t >T_a \wedge T_c)},
\end{split}
\end{equation}
and to complete the proof of  Proposition \ref{mart} we  will show that  
\begin{equation}
\label{conditional}
{\mathbb E}^b \bigl[ F| {\cal L}\bigr] = G(b) \qquad  \text{ a.s.},
\end{equation}
 for all $b \in (a,c)$, where $G=G(b)$ is given by
\begin{equation}
G=  \exp
\left\{\frac{-\kappa^2}{2}\int_a^b \frac{dy}{l(T_0,y)} - \frac{\kappa(\kappa-2)}{2}\int_{b}^c \frac{dy}{l(T_0,y)} \right\}.
\end{equation}

 As a first step in establishing \eqref{conditional}, let us verify the integrated version ${\mathbb E}^b\bigl[ F\bigr] ={\mathbb E}^b\bigl[G\bigr]$.
According to the Ray-Knight theorems the local time process $\bigl(l(T_0,y) ; y \geq 0 \bigr)$ is distributed, under ${\mathbb P}^b$, as a  inhomogeneous diffusion: starting from zero it behaves as a squared Bessel process of dimension $2$ for  $y \in [0,b]$ and then behaves as a squared Bessel process of dimension $0$ for $y \in [b,\infty)$. Conditioning on $l(T_0,a)=u, l(T_0,b)=v$ and $l(T_0,c)=w$, we have  $\bigl(l(T_0,y) ; a \leq y \leq b \bigr)$ and  $\bigl(l(T_0,y) ; b \leq y \leq c \bigr)$ are independent and distributed as squared Bessel bridges of dimension $2$ and dimension $0$ respectively.  
For any dimension $d \geq 0$ we will denote the  density of
the transition density of the
$\delta$-dimensional squared Bessel process by
$q^\delta_t(x,y)$, for $x,y,t>0$. Note that if $\delta=0$
then the transition kernel also has an atom at zero of  size $q^0_t(x,0)$. It is an 
consequence of the absolute continuity relations between the squared
Bessel processes that:
\[
{\mathbb E} \left[ \exp\left\{-\frac{\kappa^2}{2}\int_0^t ds/Z_s
\right\}\Big \vert Z_0=u, Z_t=v  \Big. \right ]=
\left(\frac{u}{v}\right)^{\kappa/2} \frac{q^{2\kappa+2}_t(u,v)}{q^2_t(u,v)}
\]
if $Z$ is a squared Bessel with dimension $2$; while
\[
{\mathbb E} \left[ \exp\left\{-\frac{\kappa(\kappa-2)}{2}\int_0^t ds/Z_s
\right\}\Big \vert Z_0=v, Z_t=w  \Big. \right
]=\left(\frac{v}{w}\right)^{\kappa/2}
\frac{q^{2\kappa}_t(v,w)}{q^{0}_t(v,w)}
\]
if $Z$ has dimension $0$. Combining these  results  gives
\begin{multline}
\label{3c}
{\mathbb E}^b \left[ \exp\left\{-\frac{\kappa^2}{2}\int_a^b \frac{dy}{l(T_0,y)}
-\frac{\kappa(\kappa-2)}{2}\int_b^c \frac{dy}{l(T_0,y)}\right\}\Big \vert l(T_0,a)=u,
l(T_0,b)=v, l(T_0,c)=w  \Big. \right]= \\
\left(\frac{u}{w}\right)^{\kappa/2}
\frac{q^{2\kappa+2}_{b-a}(u,v)q^{\kappa}_{c-b}(v,w)}
{q^{2}_{b-a}(u,v)q^{0}_{c-b}(v,w)}
.
\end{multline}
 The joint law of
$l(T_0,a)$, $l(T_0,b)$ and $l(T_0,c)$ (restricted to $w>0$)  is
\[
q_a^2(0,u)q^2_{b-a}(u,v)q^{0}_{c-b}(v,w)\;du\,dv\,dw. 
\]
So integrating out, noting that there is no contribution to the expectation from the event $\{l(T_0,c)=0\}$,
\begin{multline}
\label{int}
{\mathbb E}^b \left[ \exp\left\{-\frac{\kappa^2}{2}\int_a^b \frac{dy}{l(T_0,y)}
-\frac{\kappa(\kappa-2)}{2}\int_b^c \frac{dy}{l(T_0,y)}\right\} \right]= \\
\int_0^\infty du \int_0^\infty dv \int _0^\infty dw \; \left(\frac{u}{w}\right)^{\kappa/2}
q_a^2(0,u)q^{2\kappa+2}_{b-a}(u,v)q^{2\kappa}_{c-b}(v,w).
\end{multline}
\begin{lemma}
For any $u,w,s,t>0$, and any $\delta \geq 0$,
\begin{equation*}
\int_0^\infty  dv\; q^{\delta+2}_s(u,v) q^\delta_t(v,w) =
\frac{s}{s+t}q^{\delta+2}_{s+t}(u,w) +\frac{t}{s+t}q^{\delta}_{s+t}(u,w).
\end{equation*}
\end{lemma}
\begin{proof}
The  identity can be rewritten as an identity-in-law. Let
$Z^\delta_t(u)$ denote a random variable distributed as the value at time
$t$ of a $\delta$-dimensional squared Bessel process starting from $u$. Then the claim is that
\[
Z^\delta_t\bigl(
Z^{\delta+2}_s(u)
\bigr)\stackrel{law}{=}1_{(I=0)}Z^{\delta+2}_{s+t}(u)+1_{(I=1)}Z^{\delta}_{s+t}(u)
\]
where $I$  is a suitable Bernoulli
 random variable independent
of other variables. Assume  this holds for $\delta=0$, then several  applications of the additive property,
\[
Z^{\delta+\delta^\prime}_t(x+y)\stackrel{law}{=} Z^\delta_t(x)+
Z^{\delta^\prime}_t(y)
\]
give us
\[
Z^\delta_t\bigl(
Z^{\delta+2}_s(u)
\bigr)\stackrel{law}{=}Z^\delta_t\bigl(
Z^{2}_s(u)+Z_s^{\delta}(0)
\bigr)\stackrel{law}{=}Z^0_t\bigl(
Z^{2}_s(u)
\bigr)+Z^\delta_{s+t}(0)\stackrel{law}{=}1_{(I=0)}Z^{\delta+2}_{s+t}(u)+1_{(I=1)}Z^{\delta}_{s+t}(u).
\]
Finally the $\delta=0$ case is easily explained by  combining the
Ray-Knight theorems with hitting probabilities for Brownian motion.
\end{proof}
Applying this lemma to the expression  at \eqref{int} gives 
\begin{equation}
\label{lin}
{\mathbb E}^b[G] =\frac{b-a}{c-a}{\mathbb E}^c\bigl[ G_+\bigr]+
\frac{c-b}{c-a} {\mathbb E}^a\bigl[ G_-\bigr],
\end{equation}
where 
\begin{align}
G_+= &\exp \left\{\frac{-\kappa^2}{2}\int_a^c \frac{dy}{l(T_0,y)}\right\}, \\
G_-=&\exp \left\{ - \frac{\kappa(\kappa-2)}{2}\int_{a}^c \frac{dy}{l(T_0,y)} \right\}.
\end{align}
On the other hand, applying the strong Markov property of Brownian motion at $T_a \wedge T_c$ and using ${\mathbb P}^b \bigl(T_a<T_c)= {c-b}/{c-a}$ we obtain \begin{equation}
{\mathbb E}^b\bigl[ F \bigr] = 
\frac{b-a}{c-a}{\mathbb E}^c\bigl[ G_+\bigr]+
\frac{c-b}{c-a} {\mathbb E}^a\bigl[ G_-\bigr].
\end{equation}
Thus, comparing this with the preceding \eqref{lin},  we have obtained
\begin{equation}
\label{equality}
{\mathbb E}^b[ F] = {\mathbb E}^b[G] \qquad \text{ for all } a\leq b \leq c.
\end{equation}

Next we introduce a family of transformations. Suppose that
$h:[0,\infty) \mapsto [0,\infty)$ is   increasing with strictly positive derivative, and  satisfies
$h(0)=0$, $h(\infty)=\infty$. If $W$ is path beginning at $W_0>0$  let $T_h W \equiv W^h$  denote the transformed path satisfying
\begin{equation}
h\bigl(W^h_u\bigr) = W_{H_u}, \qquad \text{ for } H_u<T_0,
\end{equation}
and $W^h_u=W_{H_u}$ for $H_u\geq T_0$, 
where $u\mapsto H_u$ is the inverse of the increasing process
\[t\mapsto \int_0^{t \wedge T_0} \frac{ds}{ (h^\prime\circ h^{-1} (W_s))^2} + (t-T_0)^+.
\]  
Let $\nu$ be a positive, finite measure supported on $(0,\infty)$. The Sturm-Liouville equation
\[
\phi^{\prime\prime}= \phi\nu
\]
admits a unique, strictly positive, decreasing solution ( see Revuz-Yor \cite{ry}, appendix)
 satisfying $\phi(0)=1$. We shall denote this solution by $\Phi_\nu$. It has the following probabilistic characterization:
 \begin{equation}
 \Phi_\nu(x)= {\mathbb E}^x\bigl[ {\cal E}_\nu (l(T_0,\cdot))\bigr] \qquad \qquad \text{ for } x\geq 0.
 \end{equation}
 Define the probability measure ${\mathbb P}^{b,\nu}$ via
\[
{\mathbb P}^{b,\nu}= \frac{1}{\Phi_\nu(b)}{\cal E}_{\nu}(l(T_0,\cdot)) \cdot {\mathbb P}^b.
\]
\begin{lemma}
\label{changeofmeasure}
Let $\nu$ be a positive, finite measure supported on $(0,\infty)$. Take $h$ to be  given by
\[
 h(y)= \int_0^y \frac{dx}{\Phi_\nu(x)^{2}}.
 \]
 Then if $W$ has law  ${\mathbb P}^{h(b)}$, the transformed process $T_hW$ has law ${\mathbb P}^{b,\nu}$.
 \end{lemma}
 \begin{proof}
Under ${\mathbb P}^b$ the density
 \[
 \frac{1}{\Phi_\nu(b)}
  \exp\left\{-{\textstyle{\frac{1}{2}}}\int_0^\infty \nu(dy) l(T_0,y)\right\} \]
  is the  terminal value  of the exponential  martingale
\begin{multline*}
  { \frac{\Phi_\nu(W_{t\wedge T_0})}{\Phi_\nu(b)}} 
  \exp\left\{ -{\textstyle{\frac{1}{2}}}\int_0^\infty \nu(dy) l(t\wedge T_0,y)\right\} = \\
  \exp\left\{\ln\Phi_\nu(W_{t \wedge T_0})-\ln\Phi_\nu(W_0)  -{\textstyle{\frac{1}{2}}}\int_0^\infty \nu(dy) l(t\wedge T_0,y)\right\} = \\
  \exp\left\{ \int_0^{t\wedge T_0}  { \frac{\Phi^\prime_\nu(W_s)}{\Phi_\nu(W_s)}} dW_s - 
  {\textstyle \frac{1}{2}}  \int_0^{t \wedge T_0}  \left( \frac{\Phi^\prime_\nu(W_s)}{\Phi_\nu(W_s)} \right)^2 ds\right\}.
  \end{multline*}
Whence, by Girsanov's formula, under ${\mathbb P}^{b, \nu}$, 
\[
W_{t}- \int_0^{t\wedge T_0}  { \frac{\Phi^\prime_\nu(W_s)}{\Phi_\nu(W_s)}} ds
\]
is a martingale. Now  we just need to notice that the scale function corresponding to this drift $b(\cdot)= \Phi^\prime_\nu(\cdot)/\Phi_\nu(\cdot)$ is  proportional to 
\[
 \int_0^y dx\; \exp\left\{ -\int_0^x 2b(z)dz \right\} =  \int_0^y \frac{dx}{\Phi_\nu(x)^{2}}.
\]
Standard scale/speed measure arguments complete the proof.
 \end{proof}

\begin{lemma}
\label{localchange}  If the occupation measure of the path $W$ admits local times  $l(t,y)$  then the transformed path $T_h W$  admits as local times 
 $l^h(\cdot,\cdot)$ given by
\[
l^h(t,y)=\frac{1}{h^\prime(y)}  l (H_t, h(y)), \qquad \text{ for } 0\leq H_t\leq T_0.
\]
 In particular $T^h_0$, the time at which the path $W^h$ first reaches zero, satisfies  $H_{T^h_0}=T_0$ and
\[
l^h(T^h_0,y)= \frac{1}{h^\prime(y)} l (T_0,h(y)).
\]
\end{lemma}
\begin{proof}
For any test function $f$, and $t $ satisfying $ H_t\leq T_0$,
\begin{multline*}
\int_0^{t} f(W^h_s)ds= \int_0^t f\circ h^{-1}\bigl(W_{H_s}\bigr) ds
= \int_0^{H_t}\frac{ f\circ h^{-1}(W_s)}{\bigl(h^\prime\circ h^{-1}(W_s)\bigr)^2} ds \\
= \int_0^\infty\frac{ f\circ h^{-1}(x)}{\bigl(h^\prime\circ h^{-1}(x)\bigr)^2} l(H_t,x)dx 
= \int_0^\infty \frac{f(y)}{h^\prime(y)} l(H_t,h(y))dy.
\end{multline*}
\end{proof}
Treating $T_h$ as an almost everywhere defined application from $\Omega$ to $\Omega$, we have with obvious notation,
\begin{multline}
F^h=F\circ T_h = \\
 {\mathbf 1}\bigl( T^h_c<T^h_a\bigr) \exp
\left\{ \frac{-\kappa^2}{2}\int_a^{c} \frac{dy}{\lambda^h(T^h_c,y)} \right\}  
+  {\mathbf 1}\bigl( T^h_a<T^h_c\bigr)  \exp\left\{-\frac{\kappa(\kappa-2)}{2}\int_{a}^c \frac{dy}{\lambda^h(T^h_a,y)} \right\}=\\
 {\mathbf 1}\bigl(( T_{h(c)}<T_{h(a)}\bigr) \exp
\left\{ \frac{-\kappa^2}{2}\int_{h(a)
}^{h(c)} \frac{dy}{\lambda(T_{h(c)},y)} \right\}  
 +  {\mathbf 1}\bigl( T_{h(a)}<T_{h(c)} \bigr) \exp\left\{-\frac{\kappa(\kappa-2)}{2}\int_{h(a)}^{h(c)} \frac{dy}{\lambda(T_{h(a)},y)} \right\}.
\end{multline}
Similarly we have
\begin{multline}
G^h=G \circ T_h =  \exp
\left\{\frac{-\kappa^2}{2}\int_a^b \frac{dy}{l^h(T^h_0,y)} - \frac{\kappa(\kappa-2)}{2}\int_{b}^c \frac{dy}{l^h(T^h_0,y)} \right\} \\
= \exp
\left\{\frac{-\kappa^2}{2}\int_{h(a)}^{h(b)} \frac{dy}{l(T_0,y)} - \frac{\kappa(\kappa-2)}{2}\int_{h(b)}^{h(c)} \frac{dy}{l(T_0,y)} \right\}.
\end{multline}
Since the arguments leading to \eqref{equality} hold  equally if we replace throughout  $a$ by $h(a)$, $b$ by $h(b)$  and $c$ by $h(c)$,   we deduce that
\begin{equation}
\label{equalityforall}
{\mathbb E}^{h(b)}\bigl[ F^h\bigr] ={\mathbb E}^{h(b)}\bigl[ G^h\bigr].
\end{equation}
Now suppose that $h$ and $\nu$ are associated as in Lemma \ref{changeofmeasure}.
Then we have, by virtue of \eqref{equalityforall}, 
\begin{equation}
\frac{1}{\Phi(b)}{\mathbb E}^b\bigl[ F {\cal E}_\nu(l(T_0, \cdot))\bigr] = 
{\mathbb E}^{b,\nu} \bigl [ F\bigr] = {\mathbb E}^{h(b)}\bigl[ F^h\bigr]= {\mathbb E}^{h(b)}\bigl[ G^h\bigr]= {\mathbb E}^{b,\nu} \bigl [ G\bigr]= \frac{1}{\Phi(b)}{\mathbb E}^b\bigl[ G {\cal E}_\nu(l(T_0, \cdot))\bigr].
\end{equation}
The measure $\nu$ being arbitrary this proves that ${\mathbb E}^b\bigl[ F| {\cal L}\bigr]=G$ and the proof of Theorem \ref{main} is complete.

\section{Ray-Knight theorems}

Suppose the Bass-Burdzy flow  $X$ is constructed from a Brownian motion $W$ as in Theorem \ref{main}. 
 Then 
 \begin{equation}
\label{convert}
 DX_{A_t}( \cdot)=\frac{\phi_t^\prime \circ \phi_0^{-1}(\cdot)}{
\phi^\prime_0 \circ \phi_0^{-1}(\cdot)} = \frac{\lambda(0, \phi_0^{-1}(\cdot))}{
\lambda(t, \phi_0^{-1}(\cdot))}.
\end{equation}
This allows us to convert  statements concerning  the distribution of the derivative of the flow into statements concerning the local times of $W$ and vice-versa. In this section we illustrate this by giving a proof of one of the Ray-Knight theorems for the flow by means of verifying the corresponding statement regarding Brownian local times.

The following result   is part of Theorem 1.1 of \cite{hw}, together with equation (1.18) there. 
Notice that, for each $x \in {\mathbb R}$, the process  $t\mapsto DX_t(x)$ is constant except when $t\mapsto X_t(x)$ visits zero.  In the case  $\beta_2=1$ and $\beta_1=0$  then $t \mapsto X_t(x)$ is transient in the sense $X_t(x)\rightarrow \infty $ as $t \rightarrow \infty$ and so $t\mapsto DX_t(x)$
is eventually constant. It is thus meaningful to consider the limit 
$DX_\infty(x)$.
\begin{theorem}
\label{rk}
Suppose that $(X_t; t\geq 0  )$ is the Bass-Burdzy flow with $\beta_1=0$ and $\beta_2=1$.  Then $1/DX_\infty(0)$ is uniformly distributed on $[0,1]$. Conditionally on $\{1/DX_\infty(0)=y\}$ the processes 
$\bigl(Y^+_x; x\geq 0\bigr)$ and $\bigl(Y^-_x; x\geq 0\bigr)$ given by
\[
Y^+_x= \frac{1}{DX_\infty(x)} \qquad \text{ and } \qquad Y^-_x=\frac{1}{DX_\infty(-x)},
\]
are  independent diffusions on $[0,1]$,  starting from $y$, with infinitesimal generators 
\[
2y(1-y) \frac{d^2}{dy^2} +2(1-y) \frac{d}{dy}\qquad \text{ and }\qquad
2y(1-y) \frac{d^2}{dy^2} +\bigl(2(1-y)-2y\bigr) \frac{d}{dy},
\]
respectively. 
\end{theorem}

It it should be noted that versions of this result hold for other values of $\beta_1$ and $\beta_2$, but by a combination of scaling and Girsanov transformations they may be deduced from the special case stated here.

The diffusions appearing in the above description of $DX_\infty$ are sometimes called Jacobi diffusions.  A Jacobi diffusion with dimensions $d_1$ and $d_2$ 
has generator
\[
2y(1-y) \frac{d^2}{dy^2} +\bigl(d_1(1-y)-d_2 y\bigr) \frac{d}{dy}
\]
The behaviour at the boundaries $0$  and $1$ is determined from the corresponding dimension, $d_1$ for the boundary at $0$ and $d_2$ for the boundary at $1$, according to the following rule. If $d=0$ the boundary is absorbing, if $0<d<2$ then the boundary is instantaneously reflecting, while  if $d\geq 2$ then the boundary point is an inaccessible entrance point.

Next we give a description of Brownian local times  corresponding to Theorem \ref{rk}.  Recall that the Brownian motion $W$ starts from $W_0=\xi>0$,  that $M= \sup\{ W_t;t\in [0, T_0]\}$, and $T_M$ is the  almost surely unique time such that $W_{T_M}=M$. Now as was noted earlier
$A_{T_M-}=\infty$. So  if we take $t=T_M$ in  equation \eqref{convert} we obtain
\begin{align}
\label{yminus}
Y^-\left(\int_\xi^x \frac{dz}{l(T_0,z)}\right) &= \frac{l(T_0,x)-l(T_M,x)}{l(T_0,x)}  \qquad \text{ for } \xi\leq x< M, \\
\label{yplus}
Y^+\left(\int_x^\xi \frac{dz}{l(T_0,z)}\right) &= \frac{l(T_0,x)-l(T_M,x)}{l(T_0,x)}  \qquad \text{ for } 0< x\leq \xi, 
\end{align}
where $Y^+$ and $Y^-$ are defined as  in Theorem \ref{rk}.   Thus Theorem \ref{rk} is equivalent to the following proposition.
\begin{proposition}
\label{localrk}
Suppose that processes $Y^-$ and $Y^+$ are determined  from the local times of a Brownian motion via equations \eqref{yminus} and \eqref{yplus}. Then $Y^-_0=Y^+_0$ is uniformly distributed on $[0,1]$ and conditionally on $Y^-_0=Y^+_0=y$ the processes $Y^-$ and $Y^+$ are independent Jacobi diffusions with dimensions $d^-_1=d^-_2=2$ and $d^+_1=2, d^+_2=0$ respectively.
\end{proposition}

The principal tool used in proving this proposition is the skew-product involving Jacobi processes and squared Bessel processes, see \cite{wy}. Suppose that  $\bigl(Z(t);t\geq 0\bigr)$ and  $\bigl(Z^\prime(t);t\geq
0\bigr)$ are independent squared Bessel processes with dimensions $d$ and
$d^\prime$ respectively and that $d+d^\prime>0$. Let
$Z^{(+)}(t)=Z(t)+Z^\prime(t)$ and suppose $Z^{(+)}(0)>0$. Define $Y$
via,
\begin{equation}
\label{skew}
 Y\left( \int_0^t ds/Z^{(+)}(s) \right) =\frac{Z(t)}{Z^{(+)}(t)}.
\end{equation}
 Then $Y$  is a Jacobi process with dimensions $d$ and
$d^\prime$ independent of $Z^{(+)}$ which is a squared Bessel process of
dimension $d+d^\prime$. In fact, as was remarked in \cite{wy}, the skew product  holds also in the case that  $Z$ and $Z^\prime$ are replaced  by $\tilde Z$ and $\tilde{Z}^\prime$ which are obtained from $Z$ and $Z^\prime$ by
\begin{equation}
\tilde{Z}_t= \frac{1}{u^\prime(t)} Z_{u(t)} \qquad \text{ and } \qquad 
\tilde{Z}^\prime_t= \frac{1}{u^\prime(t)} Z^\prime_{u(t)}
\end{equation}
where $u$ is a strictly increasing, continuously differentiable function satisfying $u(0)=0$. 
We  will make use of this in the case $u(t)=t/(h-t)$, in which case 
$ \tilde{Z}$ and $\tilde{Z^\prime}$ are independent squared Bessel bridges leading to $0$  over the interval $[0,h]$. 

We may now proceed with the proof of Proposition \ref{localrk}. We  decompose the path $(W_t ;0\leq t \leq T_0)$ about its maximum: let
\begin{align}
W^{(1)}_t &= W_t   &\text{ for } 0\leq t\leq T_M \\
W^{(2)}_t &= W_{T_0-t}  &\text{ for } 0\leq t \leq T_0-T_M.
\end{align}
Then conditionally on $M=m\in (\xi,\infty)$ the path segments $(W^{(1)}_t; 0 \leq t\leq T_M)$ and $(W^{(2)}_t; 0 \leq t\leq T_0-T_M)$ are independent Bessel processes of dimension  three  run until  first hitting the level $m$, started from $\xi$ and $0$ respectively. Using this decomposition we obtain the description of the local times of $W$ in the following lemma, and then
Proposition \ref{localrk}, and hence Theorem \ref{rk}, follow by virtue of the skew product.

 We  use the notation ${\mathbb Q}^{d,h}_{x,y}$ to denote the law of a squared Bessel bridge of dimension $d$ leading from $x$ to $y$ over the interval $[0,h]$.
\begin{lemma}
\label{classic}
 Conditionally on $M=m$, the local time processes $l(T_M, \cdot)$ and $l(T_0, \cdot)- l(T_M, \cdot)$ are distributed as follows.
\begin{description}
\item $\bigl(l(T_M,x); 0\leq x \leq m\bigr)$ and 
$\bigl(l(T_0,x)-l(T_M,x); 0\leq x \leq m\bigr)$ are independent.
\item $l(T_M, \xi)$ and $l(T_0,\xi)-l(T_M,\xi)$ each have the exponential distribution with mean $2\xi(m-\xi)/m$. 
\item Conditionally on $l(T_M, \xi)= z$ the processes 
$\bigl( l(T_M, \xi-x); 0\leq x\leq \xi\bigr)$ and $\bigl( l(T_M, \xi+x); 0\leq x\leq m-\xi\bigr)$ are independent   and are distributed as ${\mathbb Q}^{0,\xi}_{z,0}$ and  ${\mathbb Q}^{2,m-\xi}_{z,0}$ respectively.
\item   Conditionally on $l(T_0, \xi)-l(T_M,\xi)= z$ the processes 
$\bigl( l(T_0,\xi-x)-l(T_M, \xi-x); 0\leq x\leq \xi\bigr)$ and $\bigl(l(T_0,\xi+x)- l(T_M, \xi+x); 0\leq x\leq m-\xi\bigr)$ are independent   and are distributed as ${\mathbb Q}^{2,\xi}_{z,0}$ and  ${\mathbb Q}^{2,m-\xi}_{z,0}$ respectively.
\end{description}
\end{lemma}
\begin{proof}
 $\bigl(l(T_M,x); 0\leq x \leq m\bigr)$ are the local times of the path segment $W^{(1)}$ while  
$\bigl(l(T_0,x)-l(T_M,x); 0\leq x \leq m\bigr)$  are the local times of $W^{(2)}$. Thus their independence follows from that of $W^{(1)}$ and $W^{(2)}$. 

Next it is well known  (see \cite{dmy}  for instance) that the local times of a
Bessel three process starting from zero taken at its hitting time of a level $m$ are distributed as ${\mathbb Q}^{2,m}_{0,0}$. So the distribution of $l(T_0,\xi)-l(T_M,\xi)$ is obtained as that of this bridge at time $\xi$:
\[
\frac{q^2_{\xi}(0,z)q^2_{m-\xi}(z,0)}{q^2_m(0,0)}= \frac{m}{2\xi(m-\xi)}\exp\left\{-\frac{zm}{2\xi(m-\xi)}\right\} ,
\]
where $q^d_t(x,y)$ are the   transition densities of a squared Bessel process of dimension $d$.
Moreover, by standard Markovian properties of bridges, a process distributed as ${\mathbb Q}^{2,m}_{0,0}$ can be constructed by conditioning on its  value at time $\xi$ being $z$, and placing ``back to back'' two independent bridges distributed as ${\mathbb Q}^{2,m-\xi}_{z,0}$ and   ${\mathbb Q}^{2,\xi}_{z,0}$. This proves that the distribution of  $\bigl(l(T_0,x)-l(T_M,x); 0\leq x \leq m\bigr)$ is as asserted.

 Turn now to the process  $\bigl(l(T_M,m-x); 0\leq x \leq m\bigr)$  which is distributed of the  local times of a Bessel three process starting from $\xi$ and  taken at its hitting time of a level $m$. This is  obtained by conditioning  the inhomogeneous diffusion which starts from $0$, evolves a squared Bessel process of dimension two until time $m-\xi$, and then evolves as a squared Bessel process of dimension zero, to have been absorbed at zero by time $m$. To see this, just note that according to the  Ray-Knight theorems, Brownian motion started from $\xi$ and taken at its hitting time of $m$ has local times distributed as the unconditioned diffusion, and that the  conditioning corresponds exactly to conditioning the path of the  Brownian motion not to have reached zero before hitting level $m$. Thus the distribution of $l(T_M,\xi)$ is given by
\[
\frac{m}{\xi}q^2_{m-\xi}(0,z)q^0_{\xi}(z,0)= \frac{m}{2\xi(m-\xi)}\exp\left\{-\frac{zm}{2\xi(m-\xi)}\right\}.
\] 
Finally, once again, by standard Markovian arguments we observe that this conditioned, inhomogeneous diffusion, can be constructed by conditioning on its value at time $\xi$ being $z$ and   placing ``back to back'' two independent bridges distributed as ${\mathbb Q}^{0,\xi}_{z,0}$ and   ${\mathbb Q}^{2,m-\xi
}_{z,0}$.
\end{proof}
\vspace{.2 in}
{\bf Acknowledgment:} Proposition \ref{mart} is a reworking of unpublished joint work with Marc Yor. This paper was written during a visit to the Department of Statistics at the University of California, Berkeley, and I would like to thank Jim Pitman and the Department for their hospitality.


\begin{thebibliography}{99}
\bibitem{aldous}
  D. Aldous, {\em Brownian excursion conditioned on its local time.}
  Elect. Comm. in Probab. {\bf 3}, 79-90, (1998).

   
\bibitem{bb}Bass, R.F. and Burdzy, K., {\em Stochastic bifurcation models},
Annals of Probability, {\bf 27}, 50-108, (1999).
   
\bibitem{dmy} 
C. Donati-Matin and M. Yor, { Some Brownian functionals and their laws.} Annals of Probability, {\bf 25}:3, 1011-1058. (1997).

\bibitem{hw}
Hu, Y. and Warren, J.,
{\em Ray-Knight theorems for a stochastic flow},
{Stochastic processes and its applications}, {\bf 86}, {287-305},
(2000).

\bibitem{j} 
Jacod, J., 
{\em  Grossissement initial, hypoth\`{e}se (H'), et th\'eor\`me de Girsanov.}
In Grossissements de filtrations: exemples et applications, 
Lecture notes in Mathematics, {\bf 1118}. Springer, (1985).
 
\bibitem{ry}
   D.Revuz and M.Yor, {\em Continuous martingales and Brownian motion},
Springer, (1998).


\bibitem{wy}
{Warren, J. and Yor, M.},
{\em The {B}rownian burglar: conditioning {B}rownian motion by its local time 
process}, {Seminaire de Probabilit\'{e}s} 32, 328-342, Lecture notes in Mathematics, {\bf 1709}, Springer, (1998).



\end{thebibliography}
\end{document}